# A REMARK ON AN INTEGRAL CHARACTERIZATION OF THE DUAL OF BV

NICOLA FUSCO

*Dipartimento di Matematica e Applicazioni Renato Caccioppoli, Università degli Studi di Napoli Federico II, Napoli, Italia*

DANIEL SPECTOR

*Department of Applied Mathematics, National Chiao Tung University, Hsinchu, Taiwan*

ABSTRACT. In this paper, we show how under the continuum hypothesis one can obtain an integral representation for elements of the topological dual of the space of functions of bounded variation in terms of Lebesgue and Kolmogorov-Burkill integrals.

## 1. INTRODUCTION AND MAIN RESULTS

Let $\Omega \subset \mathbb{R}^d$ be open and denote by $BV(\Omega)$ the set of functions of bounded variation in $\Omega$, that is, those functions $u \in L^1(\Omega)$ with finite total variation:

(1.1) $$|Du|(\Omega) := \sup_{\substack{\Phi \in C_c^1(\Omega;\mathbb{R}^d),\\ \|\Phi\|_{C_0(\Omega;\mathbb{R}^d)} \leq 1}} \int_\Omega u \operatorname{div} \Phi < +\infty.$$

When equipped with the norm
$$\|u\|_{BV(\Omega)} := \|u\|_{L^1(\Omega)} + |Du|(\Omega),$$

$BV(\Omega)$ is a Banach space, though in contrast to the Sobolev case, smooth compactly supported functions fail to be dense in the norm topology, nor is the space even separable. Its importance then lies in the fact that it is in some sense the natural space in a variety of problems in the calculus of variations where energies exhibit linear bounds in the gradient of the field - a notable example is in the study of minimal surfaces (see [4] for the systematic study of the space in this context).

A key tool in its study is the compactness of bounded sets with respect to the weak-star convergence of measures provided by the Banach-Alaoglu theorem, as $BV(\Omega)$ can be embedded in the product space of integrable functions crossed with $\mathbb{R}^d$-valued finite Radon measures $L^1(\Omega) \times M_b(\Omega;\mathbb{R}^d)$. More than this, $BV(\Omega)$ is itself a dual space (a result we describe in Section 2 in more detail) for which one has an integral representation of the duality pairing. However, a problem which has not been clearly explained is that of giving an integral characterization of its dual. We were curious to see if this could be done, as the literature on the space itself has only partial known results. The problem can be dated to at least the 1977 paper of Meyers and Ziemer [8] on Poincaré-Wirtinger inequalities that provides a characterization of the positive measures in the dual of $BV(\Omega)$, while such a question is raised explicitly in the 1984 AMS summer meeting on Geometric Measure

Theory and the Calculus of Variations (see [2], Problem 7.4 on p.458). It is again stated in the 1998 paper of De Pauw [3], who gives an integral representation for the dual of the space of special functions of bounded variation. More recently, the problem has been discussed in the papers of Phuc and Torres [9] and Torres [10].

It turns out that one can give an integral representation of the dual of $BV(\Omega)$, and that the solution existed more or less before any reference to such a problem. Indeed, in the papers [6,7], assuming the continuum hypothesis, Mauldin characterized the dual of spaces of finite Radon measures in the case the set of all Radon measures on the space has cardinality at most $2^{\aleph_0}$. As a consequence of his result and the Hahn-Banach theorem, one can easily get the following integral representation for the dual of the $BV(\Omega)$. To this end, let us recall that given a function $u \in BV(\Omega)$, the measure derivative $Du$ can be written as

$$Du = \nabla u \mathcal{L}^d + D^j u + D^c u,$$

where $\nabla u$ is absolutely continuous with respect to the Lebesgue measure $\mathcal{L}^d$, $D^j u = \nu_u(u^+ - u^-)\mathcal{H}^{d-1} \llcorner J_u$ with $J_u$ a $(d-1)$-rectifiable set, and $D^c u = D^s u - D^j u$ for $D^s u$ the portion of $Du$ which is singular with respect to the Lebesgue measure.

**Theorem 1.1.** *Let $\Omega \subset \mathbb{R}^d$ be open and $L \in (BV(\Omega))'$. Under ZFC set theory and the Continuum Hypothesis, there exists $g \in L^\infty(\Omega)$, $G_0 \in L^\infty(\Omega;\mathbb{R}^d)$, an $\mathcal{H}^{d-1}$ bounded and measurable function $G_1 : \Omega \to \mathbb{R}^d$ and $\Psi : \mathcal{B}(\Omega) \to \mathbb{R}$ a bounded, Borel set function for which*

$$L(u) = \int_\Omega gu\, dx + \int_\Omega G_0 \cdot \nabla u\, dx + \int_{J_u} G_1 \cdot dD^j u + (K)\int \Psi \cdot dD^c u$$

*for all $u \in BV(\Omega)$. Here, the last integral on the right hand side is understood in the Kolmogorov-Burkill sense. Conversely, any such integral functional is in the dual of $BV(\Omega)$.*

Here we observe that the integral representation is in terms of both Lebesgue and Kolmogorov-Burkill integrals. We recall that given a set function $\psi : \mathcal{B}(\Omega) \to \mathbb{R}$ and a countably additive set function $\mu : \mathcal{B}(\Omega) \to [-\infty, \infty]$, the Kolmogorov-Burkill integral is defined as follows. We say that $\psi$ is Kolmogorov-Burkill integrable with respect to $\mu$ if there exists a number $I \in \mathbb{R}$ such that for every $\varepsilon > 0$ there is a finite partition of $\Omega$ into Borel sets $D$ such that if $D'$ is any refinement of $D$,

$$\left|\sum_{B \in D'} \psi(B)\mu(B) - I\right| < \varepsilon.$$

As discussed by Kolmogorov in [5], this notion of integration encompasses both Lebesgue and Riemann integration. Indeed, the Kolmogorov integral generalizes the idea of Leibniz on areas as set functions, a sort of Riemann integral for the set function $\psi \cdot \mu$. Notice that while $\mu$ is countably additive, $\psi$ satisfies no such property. Moreover, for any fixed $\mu \in M_b(\Omega)$, the Kolmogorov-Burkill integral projects onto a Lebesgue integral. In fact, if $\psi : \mathcal{B}(\Omega) \to \mathbb{R}$ is Kolmogorov-Burkill integral with respect to $\mu$, then there exists $g_\mu \in L^\infty(\Omega, \mathcal{B}(\Omega), |\mu|)$ such that for every $f \in L^1(\Omega, \mathcal{B}(\Omega), |\mu|)$

$$(K)\int \psi\, d(f\mu) = \int_\Omega f g_\mu\, d\mu.$$

More generally, such a representation formula holds whenever $(\Omega, \mathcal{B}(\Omega), \mu)$ is an almost decomposable measure space. This is an immediate consequence of a result by De Pauw (see Lemma 3.4 in [3]).

Although we are not able to prove that the continuum hypothesis is necessary to give the integral representation in Theorem 1.1, Theorem 3.14 in [3] suggests that it might be so. However, it is interesting to observe that in one dimension one can characterize the dual of $SBV(a,b)$ without the continuum hypothesis. Indeed, given $L \in (SBV(a,b))'$, we can separate its action on $u, u' \in L^1(\Omega)$ and $D^j u = \sum_{x \in J_u} u(x^+) - u(x^-)$. This decoupling implies that there exists two bounded Borel functions $g, g_0$ and a function $g_1 : (a,b) \to \mathbb{R}$ such that $\sup_{x \in (a,b)} |g_1(x)| < +\infty$, for which

$$L(u) = \int_\Omega gu \, dx + \int_\Omega g_0 u' \, dx + \sum_{x \in J_u} g_1(x)(u(x^+) - u(x^-)),$$

for all $u \in SBV(a,b)$.

## 2. Proofs of the Main Results

In the introduction we mention that $BV(\Omega)$ is a dual space. We here provide a proof for the convenience of the reader. If one defines the space of distributions

$$V(\Omega) := \{T \in \mathcal{D}'(\Omega) : T = \varphi_0 + \sum_{j=1}^d \frac{\partial}{\partial x_j} \varphi_j : \varphi_j \in C_0(\Omega)\}$$

equipped with the norm

$$\|g\|_{V(\Omega)} := \|\varphi_0\|_\infty + \sum_{j=1}^d \|\varphi_j\|_\infty,$$

then we show

**Proposition 2.1.** *Let $L \in (V(\Omega))'$. Then there exists $u_L \in BV(\Omega)$ such that*

$$L(T) = \int_\Omega u_L \varphi_0 \, dx + \int_\Omega \Phi \cdot dDu_L$$

*for every $T \in V(\Omega)$ such that*

$$T = \varphi_0 + \operatorname{div} \Phi.$$

*Proof.* We recall that $L \in (V(\Omega))'$ means

$$L(\alpha T + \beta S) = \alpha L(T) + \beta L(S),$$
$$|L(T)| \leq C\|T\|_{V(\Omega)},$$

for all $\alpha, \beta \in \mathbb{R}$, $S, T \in V(\Omega)$. Let $L \in (V(\Omega))'$ be given. Then taking $T = \varphi$ where $\varphi \in C_0(\Omega)$ we have

$$|L(\varphi)| \leq C\|\varphi\|_\infty,$$

so that by the Riesz representation theorem we have $L \in M_b(\Omega)$ and

$$L(\varphi) = \int_\Omega \varphi \, d\mu_L.$$

Then for every $\Phi \in C_c^1(\Omega; \mathbb{R}^d)$, we have by continuity and this representation

$$|L(\operatorname{div} \Phi)| \leq C\|\Phi\|_\infty$$

Thus, $\mu_L \in M_b(\Omega)$ is such that its distributional derivatives $\frac{\partial \mu_L}{\partial x_j} \in M_b(\Omega)$, so that $\mu_L = u_L \in BV(\Omega)$. □

In order to prove Theorem 1.1, we require a result of Mauldin. The next statement is a special case of Theorem 10 in [7].

**Theorem 2.2.** *Let $\Sigma$ be a $\sigma$-algebra in $\Omega \subset \mathbb{R}^d$ such that the set of all countably additive real-valued measures on $\Sigma$ has cardinality at most $2^{\aleph_0}$. Assume the continuum hypothesis, $2^{\aleph_0} = \aleph_1$. Then for every $L \in (M_b(\Omega, \Sigma))'$ there exists a bounded set function $\psi : \Sigma \to \mathbb{R}$ such that*

$$L(\mu) = (K) \int \psi \, d\mu$$

*for every $\mu \in M_b(\Omega, \Sigma)$.*

With this theorem we can now give the proof of Theorem 1.1.

*Proof.* First let us observe that we may identify $BV(\Omega)$ with a closed subspace of $L^1(\Omega) \times M_b(\Omega; \mathbb{R}^d)$. Indeed, the space of distributional gradients is curl-free, in the sense of distributions, and so is closed with respect to sequential convergence in $L^1(\Omega)$, and so in particular it is closed with respect to the norm topology on $BV(\Omega)$. Thus, for any $L \in (BV(\Omega))'$, we may by the Hahn-Banach theorem extend $L$ to a functional $\overline{L} \in (L^1(\Omega) \times M_b(\Omega; \mathbb{R}^d))'$ such that

$$L(u) = \overline{L}(u, Du)$$

for all $u \in BV(\Omega)$. Now, for any $\overline{L} \in (L^1(\Omega) \times M_b(\Omega; \mathbb{R}^d))'$, we may separate its action on $L^1(\Omega)$ and $M_b(\Omega; \mathbb{R}^d)$. Thus, by the classical Riesz representation theorem we have

$$\overline{L}(f, \boldsymbol{\mu}) = \int_\Omega fg \, dx + \tilde{L}(\boldsymbol{\mu})$$

for some $g \in L^\infty(\Omega)$, $\tilde{L} \in (M_b(\Omega; \mathbb{R}^d))'$, and all $f \in L^1(\Omega)$ and $\boldsymbol{\mu} \in M_b(\Omega; \mathbb{R}^d)$. We can now invoke Theorem 2.2, since the cardinality of the space of countably additive real-valued measures on $\mathcal{B}(\Omega)$ is $2^{\aleph_0}$.

To see this, it is enough to show that a finite non-negative Borel measure is determined by its values on open balls $B(x, r) \subset \Omega$ where $x \in \mathbb{Q}^d$ and $r \in \mathbb{Q}^+$. Indeed, if two finite non-negative Borel measures $\mu, \nu$ coincide on such balls, then by approximation they also coincide on any closed ball $\overline{B}(x, r) \subset \Omega$ for $x \in \Omega$ and $r > 0$. It now suffices to show that $\mu$ and $\nu$ coincide on all open subsets of $\Omega$, since by approximation they will coincide on any Borel subset of $\Omega$. Thus, let $U \subset \Omega$ be open. An application of the Vitali-Besicovitch covering theorem (see, for instance, Theorem 2.19 in [1]) to the measure $\mu + \nu$ yields a sequence of disjoint closed balls $\overline{B}(x_n, r_n) \subset U$ such that

$$(\mu + \nu)\left(U \setminus \cup_n \overline{B}(x_n, r_n)\right) = 0.$$

Thus,

$$\mu\left(U \setminus \cup_n \overline{B}(x_n, r_n)\right) = \nu\left(U \setminus \cup_n \overline{B}(x_n, r_n)\right) = 0,$$

so that

$$\mu(U) = \mu\left(\cup_n \overline{B}(x_n, r_n)\right) = \nu\left(\cup_n \overline{B}(x_n, r_n)\right) = \nu(U).$$

Thus there exist bounded Borel functions $\psi_j : \mathcal{B}(\Omega) \to \mathbb{R}$, $j = 1 \ldots d$ such that

$$\tilde{L}(\boldsymbol{\mu}) = \sum_{j=1}^d (K) \int \psi_j \, d\mu_j$$

for any $\boldsymbol{\mu} \in M_b(\Omega; \mathbb{R}^d)$, where we have $\boldsymbol{\mu} = (\mu_1, \mu_2, \ldots, \mu_d)$. Define $\Psi := (\psi_1, \psi_2, \ldots, \psi_d)$, so that we may write the preceding more compactly as

$$\overline{L}(f, \boldsymbol{\mu}) = \int_\Omega fg \, dx + (K) \int \Psi \cdot d\boldsymbol{\mu}.$$

Now observe that the restriction of $\tilde{L} : L^1(\Omega; \mathbb{R}^d) \to \mathbb{R}$ is linear and continuous, and therefore again by the classical Riesz representation theorem there exists $G_0 \in L^\infty(\Omega; \mathbb{R}^d)$ such that

$$(K) \int \Psi \cdot d(F\mathcal{L}^d) = \int_\Omega G_0 \cdot F \, dx$$

for all $F \in L^1(\Omega; \mathbb{R}^d)$. We can treat in a similar manner the restriction $\tilde{L} : L^1(\Omega, \mathcal{B}(\Omega), \mathcal{H}^{d-1}; \mathbb{R}^d) \to \mathbb{R}$, since $(\Omega, \mathcal{B}(\Omega), \mathcal{H}^{d-1})$ is an almost decomposable measure space, and so by Lemma 3.4 in [3] the map

$$\mathcal{Y} : L^\infty(\Omega, \mathcal{B}(\Omega), \mathcal{H}^{d-1}) \to (L^1(\Omega, \mathcal{B}(\Omega), \mathcal{H}^{d-1}))'$$

given by

$$\mathcal{Y}(g)(f) = \int_\Omega gf \, d\mathcal{H}^{d-1}$$

is onto. Therefore there exists an $\mathcal{H}^{d-1}$ bounded and measurable function $G_1 : \Omega \to \mathbb{R}^d$ such that

$$(K) \int \Psi \cdot d(F\mathcal{H}^{d-1}) = \int_\Omega G_1 \cdot F \, d\mathcal{H}^{d-1}.$$

for all $F \in L^1(\Omega, \mathcal{B}(\Omega), \mathcal{H}^{d-1})$.

Then restricting the representation to $u \in BV(\Omega)$, along with the decomposition

$$Du = \nabla u \mathcal{L}^d + D^j u + D^c u,$$

and linearity of the Kolmogorov-Burkill integral with respect to the measure, we finally deduce that

$$L(u) = \int_\Omega gu \, dx + \int_\Omega G_0 \cdot \nabla u \, dx + \int_\Omega G_1 \cdot dD^j u + (K) \int \Psi \cdot dD^c u.$$

□

## ACKNOWLEDGEMENTS

The authors would like to thank Cindy Chen for her assistance in obtaining some of the articles utilized in the researching of this paper. N.F. supported by PRIN 2015PA5MP7 'Calcolo delle Variazioni', D.S. supported by the Taiwan Ministry of Science and Technology under research grant 105-2115-M-009-004-MY2. Part of this work was written while N.F. was visiting NCTU with support from the Taiwan Ministry of Science and Technology through the Mathematics Research Promotion Center.

## REFERENCES

[1] L. Ambrosio, N. Fusco, and D. Pallara, *Functions of bounded variation and free discontinuity problems*, Oxford Mathematical Monographs, The Clarendon Press, Oxford University Press, New York, 2000. MR1857292

[2] J. E. Brothers, *Some open problems in geometric measure theory and its applications suggested by participants of the 1984 ams summer institute*, Proc. Symp. Pure Math. **44** (1984), 441–464. Experimental mathematics: computational issues in nonlinear science (Los Alamos, NM, 1991).

[3] T. De Pauw, *On SBV dual*, Indiana Univ. Math. J. **47** (1998), no. 1, 99–121. MR1631541

[4] E. a. Giusti, *Minimal surfaces and functions of bounded variation*, Department of Pure Mathematics, Australian National University, Canberra, 1977. With notes by Graham H. Williams, Notes on Pure Mathematics, 10. MR0638362

[5] A. N. Kolmogorov and and, *Selected works of A. N. Kolmogorov. Vol. I*, Mathematics and its Applications (Soviet Series), vol. 25, Kluwer Academic Publishers Group, Dordrecht, 1991. Mathematics and mechanics, With commentaries by V. I. Arnol′d, V. A. Skvortsov, P. L. Ul′yanov et al, Translated from the Russian original by V. M. Volosov, Edited and with a preface, foreword and brief biography by V. M. Tikhomirov. MR1175399

[6] R. D. Mauldin, *A representation theorem for the second dual of $C[0, 1]$*, Studia Math. **46** (1973), 197–200. MR0346506
[7] ______ , *The continuum hypothesis, integration and duals of measure spaces*, Illinois J. Math. **19** (1975), 33–40. MR0377008
[8] N. G. Meyers and W. P. Ziemer, *Integral inequalities of Poincaré and Wirtinger type for BV functions*, Amer. J. Math. **99** (1977), no. 6, 1345–1360. MR0507433
[9] N. C. Phuc and M. Torres, *Characterization of signed measures in the dual of BV and related isometric isomorphisms*, Annali della Scuola Normale Superiore di Pisa **XVII** (2017), no. 2, arXiv:1503.06208.
[10] M. Torres, *On the dual of BV*, Contemporary Mathematics (to appear).